\documentclass[]{amsart}
\usepackage{amsmath}
\usepackage{amssymb}
\usepackage[dvips]{graphicx}
\usepackage[colorlinks=true, linkcolor=black, citecolor=black]{hyperref}
\usepackage{pstricks}
\usepackage{verbatim}
\usepackage{cmib5c}
\usepackage{graphicx} 
\usepackage{booktabs} 
\usepackage{tikz}
\usetikzlibrary{arrows}
\tikzstyle{block}=[draw opacity=0.7,line width=1.4cm]
\usetikzlibrary{calc}
\usepackage{tkz-base,tkz-euclide}
\usepackage[active]{srcltx}
\usepackage{enumerate}
\usepackage{graphicx, float}

\begin{document}

\title{On a useful lemma that relates quasi-nonexpansive and demicontractive mappings in Hilbert spaces}

\author{Vasile Berinde}

\address{Department of Mathematics and Computer
Science\newline\indent North University Centre at Baia Mare
\newline\indent Technical University of Cluj-Napoca
\newline\indent Victoriei 76, 430072 Baia Mare Romania}
\email{vasile.berinde@mi.utcluj.ro}

\address{Academy of Romanian Scientists}
\email{vasile.berinde@gmail.com}

\setcounter{page}{1}\coordinates{33}{2024}{1}{?? - ??}
\date{16.09.2023}\undate{22.11.2023}\adate{28.12.2023}


\begin{unabstract}
We give a brief account on a basic result (Lemma \ref{lem2}) which is a very useful tool in proving various convergence theorems in the framework of the iterative approximation of fixed points of demicontractive mappings in Hilbert spaces. This Lemma relates the class of quasi-nonexpansive mappings, by one hand, and the class of $k$-demicontractive mappings (quasi $k$-strict pseudocontractions), on the other hand and essentially states that the class of demicontractive mappings, which strictly includes the class of quasi-nonexpansive mappings, can be embedded in the later by means of an averaged perturbation. From the point of view of the fixed point problem, this means that any convergence result for Krasnoselskij-Mann iterative algorithms in the class of $k$-demicontractive mappings can be derived from its corresponding counterpart  from quasi-nonexpansive mappings.
\end{unabstract}

\maketitle

\section{Introduction}
Nonexpansive  type operators are extremely important in the metric fixed point theory, both from the theoretical point of view  and especially for their large areas of applications, see \cite{Ber23a} for a very recent survey. In this note we shall refer mainly to the following classes of mappings: nonexpansive, quasi-nonexpansive, $k$-strictly pseudocontractive (in the sense of Browder and Petryshyn) and quasi $k$-strictly pseudocontractive (commonly called demicontractive), which, although largely well known, are defined in the following for the sake of completeness.
\medskip

Let $H$ be a real Hilbert space with norm and inner product denoted as usually by $\|\cdot\|$ and $\langle \cdot,\cdot\rangle$, respectively. Let $C\subset H$ be a closed and convex set and $T:C\rightarrow C$ be a self mapping. Denote by
$$
Fix\,(T)=\{x\in C: Tx=x\}
$$
the set of fixed points of  $T$. 

\begin{definition}\label{def1}
The mapping $T$ is said to be:
\smallskip

1) {\it nonexpansive} if 
\begin{equation}\label{ne}
\|Tx-Ty\|\leq \|x-y\|,\,\textnormal{ for all } x,y\in C.
\end{equation}

2) {\it quasi-nonexpansive} if $Fix\,(T)\neq \emptyset$ and 
\begin{equation}\label{qne}
\|Tx-y\|\leq \|x-y\|,\, \textnormal{ for all } x\in C \textnormal{ and } y\in Fix\,(T).
\end{equation}

3) {\it $k$-strictly pseudocontractive} of the Browder-Petryshyn type (\cite{BroP66}) if there exists $k<1$ such that
\begin{equation}\label{strict}
\|Tx-Ty\|^2\leq \|x-y\|^2+k\|x-y-Tx+Ty\|^2, \forall x,y\in C.
\end{equation}

4) {\it $k$-demicontractive} (\cite{Hicks}) or {\it quasi $k$-strictly pseudocontractive}  (see \cite{BPR23}) if $Fix\,(T)\neq \emptyset$ and there exists  a positive number $k<1$ such that 
\begin{equation}\label{demi}
\|Tx-y\|^2\leq \|x-y\|^2+k\|x-Tx\|^2,
\end{equation}
for all $x\in C$ and $y\in Fix\,(T)$. 
\end{definition}
\medskip

It is known, see the remarks following Definition \ref{def2}, that Definition \ref{def1} 4) is equivalent, in the setting of a Hilbert space, with Definition \ref{def2}, that is, \eqref{demi} is equivalent to \eqref{eq1}.

Let us denote by $\mathcal{NE}$, $\mathcal{QNE}$, $\mathcal{SPC}$ and $\mathcal{DC}$ the classes of nonexpansive, quasi-nonexpansive, $k$-strictly pseudocontractive (in the sense of Browder and Petryshyn) and quasi $k$-strictly pseudocontractive (demicontractive), respectively.

In Metrical Fixed Point Theory there was a long standing and there still exists a  steadily increasing interest for studying the existence and approximation of fixed points of mappings in all of the above four classes of mappings and in many related ones like asymptotically nonexpansive, firmly nonexpansive etc. 

Most of the literature is devoted to mappings in the classes $\mathcal{NE}$, $\mathcal{QNE}$, $\mathcal{SPC}$ but, starting with the year 2008, there was also an increasing interest for studying the mappings in the class $\mathcal{DC}$, see the very recent survey \cite{Ber23a} and especially the consistent list of references therein, of which most are also included here, for the sake of completeness, see   \cite{Abkar17}-\cite{Banta}, \cite{Ber23}-\cite{Boikanyo}, \cite{Calderon}-\cite{Pany23}, \cite{Qin}-\cite{Zong}.

In order to establish convergence theorems for fixed point iteration schemes, some authors (\cite{Mar77}, \cite{Hicks}, \cite{Moudafi},...) have used implicitly or explicitly (\cite{MarP08}, \cite{Ber23}) a lemma that relates the classes $\mathcal{QNE}$ and $\mathcal{DC}$. 

The aim of this note is to review some of the most important moments in the process of discovering and use of this Lemma in order to prove convergence theorems in the class of demicontractive operators.

\section{The complete inclusion diagram of the classes $\mathcal{NE}$, $\mathcal{QNE}$, $\mathcal{SPC}$ and $\mathcal{DC}$}

To our best knowledge, there is no any paper that includes together a diagram of the four classes of nonexpansive type mappings $\mathcal{NE}$, $\mathcal{QNE}$, $\mathcal{SPC}$ and $\mathcal{DC}$, which should clearly show by appropriate examples the complete map of the relationships existing between all of them. 

So, we are doing this is the present section, mainly for its use in this note but also for the importance itself of such a diagram.

The next two simple examples show that $\mathcal{NE}$ and $\mathcal{QNE}$ are independent sets, i.e., $\mathcal{NE}\cap\mathcal{QNE}\neq \emptyset$,  $\mathcal{NE}$ is not included in $\mathcal{QNE}$ and $\mathcal{QNE}$ is not included in $\mathcal{NE}$.

\begin{example} \label{ex0}
Let $H$ be the real line with the usual norm, $C=[0,1]$ and $T_1x=1+x,\,x\in [0,1]$. Then: 1) $T_1\in \mathcal{NE}$; 2) $Fix\,(T_1)=\emptyset$; 3) $T_1\notin \mathcal{QNE}$.
\end{example}
\begin{example} \label{ex1}
Let $H$ be the real line with the usual norm,  $C=[0,2]$ and $T_2x=2-x,\,x\in [0,2]$. Then: 1) $T_2\in \mathcal{NE}$; 2) $T_2\in \mathcal{QNE}$; 3) $Fix\,(T_2)=\{1\}$.
\end{example}

The following lemma follows immediately from Definition \ref{def1}.
\begin{lemma}\label{lem1}

\begin{equation}\label{1}
\mathcal{NE}\subseteq \mathcal{SPC};
\end{equation}

\begin{equation}\label{2}
\mathcal{QNE}\subseteq  \mathcal{DC}.
\end{equation}

\end{lemma}

By means of the next example we show that inclusion \eqref{1} is strict, i.e., $\mathcal{NE}\subsetneq \mathcal{SPC}$.
\begin{example} \label{ex2}
Let $H$ be the real line with the usual norm, $C=\left[\dfrac{1}{2}, 2\right]$ and  $T_3:C\rightarrow C$ defined by $T_3(x)=\dfrac{1}{x}, \forall x\in C$. \quad Then:  1) $Fix\,(T_3)\neq \emptyset$; \quad 2) $T_3\in \mathcal{SPC}$; \quad 3) $T_3\notin \mathcal{NE}$.
\end{example}
\begin{proof}

1) $Fix\,(T_3)= \{1\}$;
\medskip

2) By \eqref{strict}, $T_3\in \mathcal{SPC}$ if there exists  $k\in(0,1)$ such that, for all $x,y\in C$,

$$\|T_3x-T_3y\|^2\leq \|x-y\|^2+k\|x-y-T_3x+T_3y\|^2,$$

 which in our case reduces to $$\left|\dfrac{1}{x}-\dfrac{1}{y}\right|^2\leq |x-y|^2+k\left|x-\dfrac{1}{x}-y+\dfrac{1}{y}\right|^2 \Leftrightarrow 1\leq x^2y^2+k(1+xy)^2, x,y\in \left[\dfrac{1}{2}, 2\right].$$ 

By denoting $t:=xy$, it follows that $t\in \left[\dfrac{1}{4}, 4\right]$ and hence we have to prove that there exists $k>0$ such that

 $\dfrac{1-t^2}{(1+t)^2}\leq k<1$, for all $t\in \left[\dfrac{1}{4}, 4\right]$. Consider the function  
$f(t):=\dfrac{1-t^2}{(1+t)^2}$, $t\in \left[\dfrac{1}{4}, 4\right]$. 

Since $f'(t)=-\dfrac{2}{(1+t)^2}<0$, it follows that  $f$ is strictly decreasing on $\left[\dfrac{1}{4}, 4\right]$, which implies 

$f(t)\leq f\left(\dfrac{1}{4}\right)=\dfrac{3}{5}$,  for all $t\in \left[\dfrac{1}{4}, 4\right]$. 

This shows that one can choose $k=\dfrac{3}{5}$ and so, $T_3$ is $\dfrac{3}{5}$-strictly pseudocontractive.

3) Assume $T_3\in \mathcal{NE}$, i.e.,  $|T_3x-T_3y|\leq |x-y|,\,\forall x,y\in C=\left[\dfrac{1}{2}, 2\right]$ and take $x=\dfrac{1}{2}$ and $y=1$ to get $|2-1|\leq \left |\dfrac{1}{2}-1\right |$ $\Leftrightarrow$ $1\leq \dfrac{1}{2}$, a contradiction.

\end{proof}

\begin{example} \label{ex3}
Let $H$ be the real line with the usual norm and $C = [0, 2]$. Define $T_4:[0,2]\rightarrow [0,2]$ by $T_4x=\dfrac{x^2+2}{x+1}$, for all $x\in [0, 2]$. Then:
1) $Fix\,(T_4)\neq \emptyset$; 2) $T_4\in \mathcal{QNE}$; 3) $T_4\notin \mathcal{NE}$; 4) $T_4\notin \mathcal{SPC}$.
\end{example}

\begin{proof}
1) $Fix\,(T_4)=\{2\}$; 

2) For $y=2$ and $x\in [0,2]$, by \eqref{qne} we have
$$
\left |T_4x-2\right|=\left |\frac{x^2+2}{x+1}-2\right|=\frac{x}{x+1}\cdot |x-2|\leq |x-2|,\,x\in [0,2],
$$
and so $T_4\in \mathcal{QNE}$.

3) Just consider $x=0$ and $y=\dfrac{1}{3}$ in \eqref{ne} to get
$$
\frac{5}{12}=\left |T_4 0-T_4 \frac{1}{3}\right| \leq \left |0- \frac{1}{3}\right|=\frac{1}{3},
$$
a contradiction since $\dfrac{5}{12}>\dfrac{1}{3}$. So, $T_4\notin \mathcal{NE}$.

4) Assume now that $T_4\in \mathcal{SPC}$ and take  $x=0$ and $y=\dfrac{1}{3}$ in \eqref{strict} to get
$$
\left(\frac{5}{12}\right)^2=\left |T_4 0-T_4 \frac{1}{3}\right|^2\leq \left |0- \frac{1}{3}\right|^2+k\left |0-T_4 0-\left(\frac{1}{3}-T_4 \frac{1}{3}\right)\right|^2=\left(\frac{1}{3}\right)^2+k\left(\frac{3}{4}\right)^2,
$$
which is a contradiction, since $\dfrac{5}{12}>\dfrac{1}{3}$ and $k>0$. Hence $T_4\notin \mathcal{SPC}$.
\end{proof}

The next example shows that the inclusion \eqref{2} is also strict, i.e., $\mathcal{QNE}\subsetneq \mathcal{DC}$.

\begin{example} \label{ex4}
Let $H$ be the real line with the usual norm and $C = [0, 1]$. Define $T_5$ on $C$ by $T_5 x = \dfrac{7}{8}$, if $0\leq x<1$ and $T_5 1 = \dfrac{1}{4}$. Then: 

1) $Fix\,(T_5)\neq \emptyset$; 2) $T_5\in \mathcal{DC}$; 3) $T_5\notin \mathcal{NE}$; 4) $T_5\notin \mathcal{QNE}$;  5) $T_5\notin \mathcal{SPC}$. 
\end{example}
\begin{proof}
1) $Fix\,(T_5)=\left\{\dfrac{7}{8}\right\}$;

2) By taking $y=\dfrac{7}{8}$ and $ x\in [0,1)$, inequality \eqref{demi} becomes: $$|T_5x-y|^2=0\leq |x-y|^2+k|x-T_5 x|^2,$$ which obviously holds,  for any $k>0$. 

It remains to check \eqref{demi} for the case $x=1$, which yields
$$
\left|\dfrac{1}{4}-\dfrac{7}{8}\right|^2\leq \left|1-\dfrac{7}{8}\right|^2+k\left|1-\dfrac{1}{4}\right|^2
$$
and which holds true for any $k\geq \dfrac{2}{3}$. 
Hence $T_5$ is $\dfrac{2}{3}$-demicontractive. 

3) To show that $T_5$ is not quasi-nonexpansive, take $x=1$ and $y=\dfrac{7}{8}$ in \eqref{qne}, to get
$
\dfrac{5}{8}\leq \dfrac{1}{8},
$
a contradiction. Hence $T_5$ is not quasi-nonexpansive.

4)  To prove that $T_5$ is not nonexpansive  take  $x=1$ and $y=\dfrac{7}{8}$ in \eqref{ne} to get the same contradiction as above.

5) Assume $T_5$ is $k$-strictly pseudocontractive, that is, there exists $k<1$ such that \eqref{strict} holds for any $x,y\in[0,1]$. By taking $x\in [0,1)$ and $y=1$ in \eqref{strict} we have
$$
\left(\frac{5}{8}\right)^2\leq (x-1)^2+k\left(x-1-\frac{5}{8}\right)^2,\,x\in [0,1),
$$
from which, by letting $x\rightarrow 1$ we obtain $1\leq k< 1$, a contradiction.

Hence $T_5$ is not strictly pseudocontractive.
\end{proof}

Based on Lemma \ref{lem1} and Examples \ref{ex0}-\ref{ex4}, we have the following complete map of the relationships between the four sets of nonexpansive type mappings in Definition \ref{def1}.
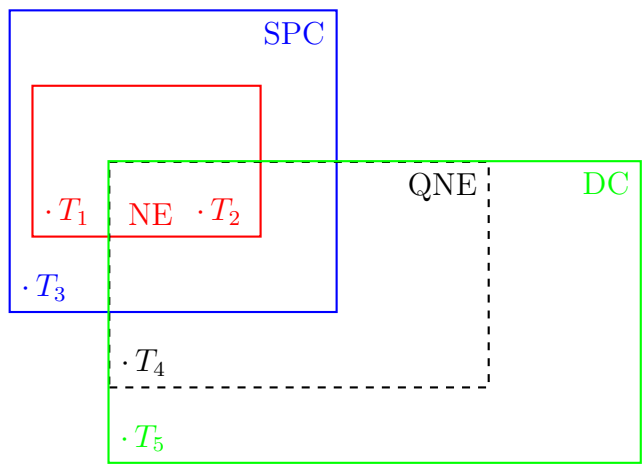
\bigskip

\begin{tikzpicture}

\draw[red, thick] node [anchor=south east]{NE} node [anchor=south west]{$\cdot \,T_2$} (-2,0) node [anchor=south west]{$\cdot \,T_1$} rectangle (1,2);
   
   \draw[blue, thick]  (-2.3,-1)  node [anchor=south west]{$\cdot \,T_3$}  rectangle (2,3) node [anchor=north east]{SPC};

      \draw[black, thick, dashed]  (-0.99,-2) node [anchor=south west]{$\cdot \,T_4$} rectangle   (4,0.99) node [anchor=north east]{QNE};
      
     \draw[green, thick]  (-1,-3) node [anchor=south west]{$\cdot \,T_5$} rectangle  (6,1)  node [anchor=north east]{DC} ;

       



\end{tikzpicture}
\medskip

{\bf Figure 1.} Diagram of the relationships between the classes $\mathcal{NE}$, $\mathcal{QNE}$, $\mathcal{SPC}$ and $\mathcal{DC}$

\bigskip

\section{A Lemma that relates quasi-nonexpansive and demicontractive mappings}

The main aim of this section is to present some historical facts about the use and formulation of an important lemma that relates quasi-nonexpansive and demicontractive mappings.

This result is of particular importance in proving convergence theorems for some fixed point iterative schemes like Krasnoselskii, Krasnoselskij-Mann etc. in the class of demicontractive mappings, by reducing the arguments to the same algorithms but in the class of quasi-nonexpansive mappings.  

We state it in the form it has been presented and used in the paper \cite{Ber23} and, for the sake of completeness, we also give its proof. 

\begin{lemma}[\cite{Ber23}, Lemma 3.2]\label{lem2}
Let $H$ be a real Hilbert space, $C\subset H$ be a closed and convex set. If $T:C\rightarrow C$ is $k$-demicontractive, then for any $\lambda\in (0,1-k)$, $T_{\lambda}$ is quasi-nonexpansive. 
\end{lemma}

\begin{proof}
By hypothesis, we have $Fix\,(T)\neq \emptyset$ and there exists $k<1$ such that 
$$
\|Tx-y\|^2\leq \|x-y\|^2+k\|x-Tx\|^2,\,x\in C \textnormal{ and } y\in Fix\,(T)
$$
which is equivalent to
$$
\langle Tx-x,x-y\rangle \leq \frac{k-1}{2}\cdot \|x-Tx\|^2,x\in C, \,y\in Fix\,(T). 
$$
Then, for all $x\in C$ and $y\in Fix\,(T)$, we have
$$
\|T_\lambda x-y\|^2=\|\lambda (Tx-x)+x-y\|^2=\|x-y\|^2+2\lambda \langle Tx-x,x-y\rangle
$$
$$
+\lambda^2 \|Tx-x\|^2\leq \|x-y\|^2+(\lambda^2+\lambda k-\lambda) \|Tx-x\|^2
$$
$$
= \|x-y\|^2+ \frac{\lambda^2+\lambda k-\lambda}{\lambda^2}\cdot \|T_\lambda x-x\|^2,\,x\in C, \,y\in Fix\,(T).
$$

So, if $\lambda^2+\lambda k-\lambda<0$, that is, $\lambda <1-k$, then the above inequality implies that
$$
\|T_\lambda x-y\|^2
\leq \|x-y\|^2,\,x\in C, \,y\in Fix\,(T),
$$
i.e., that $T_\lambda$ is quasi-nonexpansive.
\end{proof}

\medskip

We are now interested to trace back on the use of this simple but important Lemma. As it has been shown in the very recent survey paper \cite{Ber23a}, the demicontractive mappings were introduced independently in 1977 by M\u aru\c ster  \cite{Mar77} and Hicks and Kubicek \cite{Hicks}, respectively, in the setting of a Hilbert space. 

The same notion has been introduced in 1973 by M\u aru\c ster \cite{Mar73}, in the particular case of $\mathbb{R}^n$, but for the case of the nonlinear equation $U(x)=0$. By simply taking $U=I-T$, one finds the same concept as the one introduced in  \cite{Mar77}. This was the reason why in the survey paper \cite{Ber23a} we have considered 1973 as the birth date of demicontractive mappings.

In order to present some facts about the early use of Lemma \ref{lem2}, we also give here M\u aru\c ster's  definition \cite{Mar77} of demicontractive mappings. It is important to note that the term "demicontractive" was coined by Hicks and Kubicek \cite{Hicks}, who introduced it by means of inequality \eqref{demi}.

\begin{definition} [M\u aru\c ster \cite{Mar77}]\label{def2}
Let  $H$ be a real Hilbert space and $C$  a closed convex subset of  $H$. A mapping $T : C \rightarrow C$ such that $Fix\,(T )\neq \emptyset$ 
    is said to satisfy {\it condition (A)}  if there exists $\lambda >0$ such that 
\begin{equation}\label{eq1}
   \langle x-Tx, x-x^*\rangle \geq \lambda \|Tx-x\|^2, \forall x\in C, x^*\in Fix\,(T). 
\end{equation}
  \end{definition} 
  
  Despite the fact that the two definitions were introduced in the same year and in very visible magazines,  it was not apparent for a rather long time  that the two inequalities \eqref{demi} and \eqref{eq1}, which involve different formulas, are actually equivalent in the setting of a Hilbert space. 

This fact was observed more than two decades later, by Moore \cite{Moore} and is based on the following identity, valid in a real Hilbert space:
 $$
 \|x-x^*\|^2+k\|x-Tx\|^2 -\|Tx-x^*\|^2 =
 2\langle x-x^*,x-Tx\rangle-(1- k)\|x-T x\|^2,  
 $$
see  \cite{Moore}) for more details.

In our recent paper \cite{Ber23}, based on Lemma \ref{lem2},  we have explicitly proven that, in Hilbert spaces, any convergence result for a Krasnoselkij type fixed point iterative algorithm in the class of demicontractive mappings can be deduced from its counterpart in the class of quasi-nonexpansive mappings.

But this fact was known and used implicitly long before by a few researchers that were working in this area. Our aim is to survey all those attempts that precede the more recent papers \cite{Tang12},  \cite{Wang13} and \cite{Ber23}, where Lemma \ref{lem2} was explicitly stated.
\medskip

1) In the proof of Theorem 1 in M\u aru\c ster \cite{Mar77}, the author used the same arguments like the ones in the proof of Lemma \ref{lem2}. 

Indeed, if we adapt the notations in \cite{Mar77} to our current ones, i.e., we denote the fixed point of $T$ by $x^*$ instead of $\xi$ and the parameter  $t_k$ involved in the Mann iteration by $t$, what M\u aru\c ster \cite{Mar77} did, see the first 4 rows on page 70, is the following 
$$
\|T_t x-x^*\|^2=\|x-x^*-t (x-Tx)\|^2=\|x-x^*\|^2-2 t \langle x-Tx,x-x^*\rangle
$$
$$
+t^2 \|x-Tx\|^2\leq \|x-x^*\|^2+t (2\lambda-t) \|x-Tx\|^2
$$
and since  $2\lambda-t>0$, it follows that
$$
\|T_t x-x^*\|\leq \|x-x^*\|,\,x\in C, \,x^*\in Fix\,(T),
$$
which means that $T_t$ is quasi-nonexpansive for $0<t<2\lambda$.

On the other hand, if we keep in mind the relationship between $\lambda$ in \eqref{eq1} and $k$ in \eqref{demi}, that is, $\lambda=\dfrac{1-k}{2}$, then we get exactly the condition on the parameter in Lemma \ref{lem2} that ensures that the averaged operator $T_t$ is quasi-nonexpansive.

As a matter of fact, in \cite{Mar77} all the above calculations were performed directly for the sequence $x_{k+1}=T_t x_k$ and not for the mapping $T_t$.

\medskip

2) In the proof of Th\' eor\' eme in \cite{Mar73}, the same arguments were used, but for the case of the nonlinear equation $U(x)=0$. By simply taking $U=I-T$, the proof actually shows that the mapping $T_\mu$ is quasi-nonexpansive for $\mu<2\eta$, where $\eta$ corresponds to $\lambda$ in \eqref{eq1}. 

Similarly to \cite{Mar77}, the author did all the calculations in \cite{Mar73} for the sequence $x_{p+1}=T_{\mu} x_p$ and not for the mapping $T_{\mu}$.
\medskip

3) In the proof of Theorem 1 in \cite{Hicks}, the authors performed similar calculations to those in \cite{Mar77} but for the sequence $v_{n+1}=T_{d_n} v_n$ and not for the averaged mapping $T_{d_n}$.
\medskip

4) In a series of papers from the period 2003-2009, see \cite{Mar03}, \cite{MarMar05}, \cite{Mar05}-\cite{MarP09}, M\u aru\c ster used Lemma \ref{lem2} and even presented a complete proof of it, but in the framework of the proof of the main result established there. For example, in \cite{Mar03}, this is done in the proof of Theorem 2.  Lemma \ref{lem2} is also explicitly stated and proved and then used to apply Theorem 1 in \cite{Mar03} (about quasi-nonexpansive mappings) to get the desired conclusion. Similar formulations of Lemma \ref{lem2} do appear under various forms in the subsequent papers by M\u aru\c ster \cite{MarMar05}, \cite{Mar05}-\cite{MarP09}.

\medskip

5) In Remark 2.1 from Moudafi \cite{Moudafi}, Lemma \ref{lem2} is explicitly stated and proven, as follows. 

"Let $T$ be a $k$-demicontractive self-mapping on $\mathcal{H}$ with $Fix\,(T)\neq \emptyset$ and set $T_w:=(1-w)I+w T$ for $w\in (0,1]$. It is obviously checked that $Fix\,(T)=Fix\,(T_w)$. Moreover, $T_w$ is quasi-nonexpansive for $w$ small enough. Indeed, given an arbitrary $(x,q)\in\mathcal{H}\times Fix\,(T)$, we have
$$
|T_w x-q|^2=|(x-q)+w(Tx-x)|^2
$$
$$
=|x-q|^2-2w\langle x-q,x-Tx\rangle+w^2 |Tx-x|^2
$$
which by (1.5) (i.e., the demicontractive condition in M\u aru\c ster's form) yields
$$
|T_w x-q|^2\leq |x-q|^2-w(1-k-w) |Tx-x|^2.
$$
Consequently, if $w\in (0,1-k]$, then $T_w$ is quasi-nonexpansive..."

This explicit statement and its proof are reproduced in Maing\' e and Moudafi \cite{Main08e} (Remark 2.1), in Maing\' e \cite{Main08a} (Remark 4.2) and 
in some other papers by the same authors.

\medskip

6) It appears that Tang et al. \cite{Tang12} were the first ones to state explicitly Lemma \ref{lem2}, by referring to Remark 2.1 from Moudafi \cite{Moudafi}.
\medskip

7) The present author, who was not aware of the implicit or explicit statements of Lemma \ref{lem2} reviewed previously, formulated it as an auxiliary result (Lemma 3.2) in \cite{Ber23}, and, based on it, presented simpler and unifying proofs for the pioneering papers by M\u aru\c ster \cite{Mar77} and Hicks and Kubicek \cite{Hicks}. 

The title of \cite{Ber23}, {\it Approximating fixed points results for demicontractive mappings could be derived from their quasi-nonexpansive counterparts}, as well as its first conclusions reproduced below should be taken into consideration by all researchers dealing  with the study of demicontractive mappings.

"1. In this paper we have shown that the convergence theorems for Mann iteration used for approximating the fixed points of demicontractive mappings in Hilbert spaces could be derived from the corresponding convergence theorems in the class of quasi-nonexpansive mappings. 

2. Our derivation is based on an imbedding technique described by Lemma 3.2,  which essentially shows that if $T$ is $k$-demicontractive, then for any $\lambda\in (0,1-k)$, $T_{\lambda}$ is quasi-nonexpansive. 

3. In this way we obtained a unifying technique of proof for various well known results in the fixed point theory of demicontractive mappings that has been illustrated  for the case of the first two classical convergence results in the class of demicontractive mappings in literature: M\u aru\c ster \cite{Mar77} and Hicks and Kubicek \cite{Hicks}."


We note that a similar technique also works for $k$-strict pseudocontractions, which can be embedded in the class of nonexpansive mappings in Hilbert spaces. This fact was first exploited by Browder and Petryshyn  \cite{BroP66}, \cite{Pet}, and also used much later by Zhou \cite{Zhou08} in the case of nonself mappings.

\section{Conclusions}

1. I this paper we gave a brief account on a basic result (Lemma \ref{lem2}) which is a very useful tool in proving various convergence theorems in the framework of the iterative approximation of fixed points of demicontractive mappings in Hilbert spaces. This lemma relates the class of quasi-nonexpansive mappings, by one hand, and the class of $k$-demicontractive mappings (or quasi $k$-strict pseudocontractions), on the other hand and essentially states that the class of demicontractive mappings, which strictly includes the class of quasi-nonexpansive mappings, can be embedded in the later by means of an averaged perturbation. 

2. From the point of view of the fixed point problem, this means that any convergence result for Krasnoselskij-Mann iterative algorithms in the class of demicontractive mappings can be derived from its corresponding counterpart established for quasi-nonexpansive mappings.

3. The nonexpansive mappings are important  in solving various problems in data science, like image recovery, machine learning, signal processing, neural networks etc. This was the reason why, in Section 2, we presented, by means of appropriate examples, the complete map of the relationships existing amongst four important such classes: nonexpansive mappings, quasi-nonexpansive mappings, strictly pseudocontractive mappings and demicontractive mappings. To our best knowledge, this is the first time such a diagram is pictured.

4. In this context, we also collected an almost complete list of references related to the study of fixed point problem in the class of demicontractive mappings, mainly taken from \cite{Ber23a}. 

5. The main message of this note for researchers working in that area is to use Lemma \ref{lem2} when dealing with convergence theorems of Krasoselskij-Man type in the class of demicontractive mappings, in order to unify and simplify the proofs.

6. One of the main aims of this note was to trace back on the awareness and use of Lemma \ref{lem2}. We thus discovered that its inception started with the pioneering works on demicontractive mappings, due to  M\u aru\c ster \cite{Mar73}, \cite{Mar77} and Hicks and Kubicek \cite{Hicks}, and that the first explicit statement and proof of this lemma in due to M\u aru\c ster \cite{Mar03}, who did it within the proof of Theorem 2 \cite{Mar03}. 

7. A similar technique works for $k$-strict pseudocontractions, which can be embedded in the class of nonexpansive mappings in Hilbert spaces, first exploited by Browder and Petryshyn  \cite{BroP66}, see also \cite{Pet}, and also used much later by Zhou \cite{Zhou08} in the case of nonself mappings, but this should be the subject of another paper.


\section*{Acknowledgements}

The research was carried out during author's short visit (December 2023) at the Department of Mathematics, King Fahd University of Petroleum and Minerals, Dhahran, Saudi Arabia. The author is grateful to Professor Monther Alfuraidan, the Chairman of Department of Mathematics, for invitation and for providing excellent facilities during his visit.


\begin{thebibliography}{10}
\providecommand{\url}[1]{\texttt{#1}}
\providecommand{\urlprefix}{URL }

\bibitem{Abbas} Abass, H. A.; Aphane, M. An algorithmic approach to solving split common fixed point problems for families of demicontractive operators in Hilbert spaces. {\it Boll. Unione Mat. Ital.} Oct 2023 (Early Access)

\bibitem{Abkar17} Abkar, A.; Shahrosvand, E. The split common fixed point problem of two infinite families of demicontractive mappings and the split common null point problem. {\it Filomat} {\bf 31} (2017), no. 12, 3859--3874.

\bibitem{Abkar17a} Abkar, A.; Shekarbaigi, M. A synthetic algorithm for families of demicontractive and nonexpansive mappings and equilibrium problems. {\it Filomat} {\bf 31} (2017), no. 19, 5891--5908.

\bibitem{Abkar15} Abkar, A.; Tavakkoli, M. A new algorithm for two finite families of demicontractive mappings and equilibrium problems. {\it Appl. Math. Comput.} {\bf 266} (2015), 491--500.

\bibitem{Adamu} Adamu, A.; Adam, A. A. Approximation of solutions of split equality fixed point problems with applications. {\it Carpathian J. Math.} {\bf 37} (2021), no. 3, 381--392.

\bibitem{Agwu} Agwu, I. K.; Igbokwe, D. I. New Iteration Algorithms for Solving Equilibrium Problems and Fixed Point Problems of Two Finite Families of Asymptotically Demicontractive Multivalued Mappings. {\it Sahand Comm. Math. Anal.} {\bf 20} (2023), No. 2, 1--38.

\bibitem{Ala} Alakoya, T. O.; Jolaoso, L. O.; Mewomo, O. T. A general iterative method for finding common fixed point of finite family of demicontractive mappings with accretive variational inequality problems in Banach spaces. {\em Nonlinear Stud.} {\bf 27} (2020), no. 1, 213--236.

\bibitem{Amarachi} Amarachi Uzor, V.; Alakoya, T. O.; Mewomo, O. T. Strong convergence of a self-adaptive inertial Tseng's extragradient method for pseudomonotone variational inequalities and fixed point problems. {\it Open Math.} {\bf 20} (2022), no. 1, 234--257.

\bibitem{Anh} Anh, T. V.; Muu, L. D.; Son, D. X. Parallel algorithms for solving a class of variational inequalities over the common fixed points set of a finite family of demicontractive mappings. {\it Numer. Funct. Anal. Optim.} {\bf 39} (2018), no. 14, 1477--1494.

\bibitem{Aremu} Aremu, K. O.; Jolaoso, L. O.; Izuchukwu, C.; Mewomo, O. T. Approximation of common solution of finite family of monotone inclusion and fixed point problems for demicontractive multivalued mappings in CAT(0) spaces. {\em Ric. Mat.} {\bf 69} (2020), no. 1, 13--34.

\bibitem{Arfat} Arfat, Y.; Kumam, P.; Khan, M. A. A.; Cho, Y. J. A hybrid steepest-descent algorithm for convex minimization over the fixed point set of multivalued mappings. {\it Carpathian J. Math.} {\bf 39} (2023), no. 1, 303--314.

\bibitem{Arfat} Arfat, Y.; Kumam, P.; Phiangsungnoen, S.; Khan, M. A. A.; Fukhar-ud-din, H. An inertially constructed projection based hybrid algorithm for fixed point and split null point problems. {\it AIMS Math.} {\bf 8} (2023), no. 3, 6590--6608.

\bibitem{Banta} Bantaojai, T.; Garodia, C.; Uddin, I.; Pakkaranang, N.; Yimmuang, P. A novel iterative approach for solving common fixed point problems in geodesic spaces with convergence analysis. {\it Carpathian J. Math.} {\bf 37} (2021), no. 2, 145--160.

\bibitem{Batra} Batra, C.; Gupta, N.; Chugh, R.; Kumar, R. Generalized viscosity extragradient algorithm for pseudomonotone equilibrium and fixed point problems for finite family of demicontractive operators. {\it J. Appl. Math. Comput.} {\bf 68} (2022), no. 6, 4195--4222.

\bibitem{Beg} Beg I.; Abbas, M.; Asghar, M. W. Approximation of the Solution of Split Equality Fixed Point Problem for Family of Multivalued Demicontractive Operators with Application. {\it Mathematics} {\bf 11} (2023), no. 4, Article number 959.

\bibitem {Ber13} Berinde, V. Convergence theorems for fixed point iterative methods defined as admissible perturbations of a nonlinear operator. {\em Carpathian J. Math.} {\bf 29} (2013), no. 1, 9--18.


\bibitem {Ber18} Berinde, V. Weak and strong convergence theorems for the Krasnoselskij iterative algorithm in the class of enriched strictly pseudocontractive operators. {\em An. Univ. Vest Timi\c s. Ser. Mat.-Inform.} {\bf 56} (2018), no. 2, 13--27.

\bibitem {Ber19a} Berinde, V. Approximating fixed points of enriched nonexpansive mappings by Krasnoselskij iteration in Hilbert spaces. {\em Carpathian J. Math.} {\bf 35} (2019), no. 3, 293--304.


\bibitem {Ber20} Berinde, V. Approximating fixed points of enriched nonexpansive mappings in Banach spaces by using a retraction-displacement condition. {\em Carpathian J. Math.} {\bf 36} (2020), no. 1, 27--34.

\bibitem {Ber23} Berinde, V. Approximating fixed points of demicontractive mappings via the quasi-nonexpansive case. {\it Carpathian J. Math.} {\bf 39} (2023), no. 1, 73--85.

\bibitem {Ber23a} Berinde, V. Single-Valued Demicontractive Mappings: Half a Century of Developments and Future Prospects. {\it Symmetry}, {\bf 15} (2023), no. 10, 1866; https://doi.org/10.3390/sym15101866


\bibitem{BerP19}  Berinde, V.; P\u acurar, M. Within the world of demicontractive mappings. In Memoriam Professor \c Stefan M\u aru\c ster (1937-2017). {\em An. Univ. Vest Timi\c s. Ser. Mat.-Inform.} {\bf 57} (2019), no. 1, 3--12.

\bibitem{BerP21c}  Berinde, V.; P\u acurar, M.  Fixed points theorems for unsaturated and saturated classes of contractive mappings in Banach spaces. {\em Symmetry} {\bf 13} (2021), Article Number 713 https://doi.org/10.3390/sym13040713.

\bibitem{BPR23} Berinde, V.; Petru\c sel, A.; Rus, I. A. Remarks on the terminology of the mappings in fixed point iterative methods in metric spaces. {\it Fixed Point Theory} {\bf 24} (2023), no. 2, 525--540.

\bibitem{Boikanyo} Boikanyo, O. A. A strongly convergent algorithm for the split common fixed point problem. {\it Appl. Math. Comput.} {\bf 265} (2015), 844--853.

\bibitem {BroP66} Browder, F. E.; Petryshyn, W. V. The solution by iteration of nonlinear functional equations in Banach spaces. {\em Bull. Amer. Math. Soc.} {\bf 72} (1966), 571--575.


\bibitem {Byrne} Byrne, C. Iterative oblique projection onto convex sets and the split feasibility problem. {\it Inverse Problems} {\bf 18} (2002), no. 2, 441--453.

\bibitem {Calderon} Calder\' on, K.; Khamsi, M. A.; Mart\' inez-Moreno, J. Perturbed approximations of fixed points of nonexpansive mappings in $\rm CAT_p(0)$ spaces. {\it Carpathian J. Math.} {\bf 37} (2021), no. 1, 65--79.

\bibitem{Chen} Chen, H. Y.; Sahu, D. R.; Wong, N. C. Iterative algorithms for solving multiple split common fixed problems in Hilbert spaces. {\it J. Nonlinear Convex Anal.} {\bf 19} (2018), no. 2, 265--285.

\bibitem{Chang} Chang, S.-S.; Wang, L.; Wang, X. R.; Zhao, L. C. Common solution for a finite family of minimization problem and fixed point problem for a pair of demicontractive mappings in Hadamard spaces. {\it Rev. R. Acad. Cienc. Exactas F\' is. Nat. Ser. A Mat. RACSAM} {\bf 114} (2020), no. 2, Paper No. 61, 12 pp.

\bibitem{Charoen} Charoensawan P.; Suparatulatorn R. A modified Mann algorithm for the general split problem of demicontractive operators.
{\it Results Nonlinear Anal.} {\bf 5} (2022), no.3 3, 213--221.

\bibitem{Chen21} Chen, H.-Y. Weak and strong convergence of inertial algorithms for solving split common fixed point problems. {\it J. Inequal. Appl.} {\bf 2021}, Paper No. 26, 17 pp.

\bibitem{Chen18} Chen, H. Y.; Sahu, D. R.; Wong, N. C. Iterative algorithms for solving multiple split common fixed problems in Hilbert spaces. {\it J. Nonlinear Convex Anal.} {\bf 19} (2018), no. 2, 265--285.

\bibitem {Chid84} Chidume, C. E. The solution by iteration of nonlinear equations in certain Banach spaces. {\em J. Nigerian Math. Soc.} {\bf 3} (1984), 57--62 (1986).

\bibitem {Chid94} Chidume, C. E. An iterative method for nonlinear demiclosed monotone-type operators. {\it Dynam. Systems Appl.} {\bf 3} (1994), no. 3, 349--355.

\bibitem {Chid09} Chidume, C. {\it Geometric properties of Banach spaces and nonlinear iterations}. Lecture Notes in Mathematics, 1965. Springer-Verlag London, Ltd., London, 2009.

\bibitem {Chid16} Chidume, C. E.; Minjibir, M. S. Krasnoselskii algorithm for fixed points of multivalued quasi-nonexpansive mappings in certain Banach spaces. {\em Fixed Point Theory} {\bf 17} (2016), no. 2, 301--311.

\bibitem{ChidM} Chidume, C. E.; M\u aru\c ster, \c S. Iterative methods for the computation of fixed points of demicontractive mappings. {\it J. Comput. Appl. Math.} {\bf 234} (2010), no. 3, 861--882.

\bibitem{ChidBN} Chidume, C. E.; Bello, A. U.; Ndambomve, P. Strong and $\Delta$-convergence theorems for common fixed points of a finite family of multivalued demicontractive mappings in ${\rm CAT}(0)$ spaces. {\it Abstr. Appl. Anal.} {\bf 2014}, Art. ID 805168, 6 pp.

\bibitem{Chid-JNAO} Chidume, C. E.; Ndambomve, P.; Bello, A. U. The split equality fixed point problem for demi-contractive mappings. {\it J. Nonlinear Anal. Optim.} {\bf 6} (2015), no. 1, 61--69.


\bibitem {Cui21} Cui, H. H. Multiple-sets split common fixed-point problems for demicontractive mappings. {\it J. Math.} {\bf 2021}, Art. ID 3962348, 6 pp.

\bibitem {Cui18} Cui, H. H.; Ceng, L. C.; Wang, F. H. Weak convergence theorems on the split common fixed point problem for demicontractive continuous mappings. {\it J. Funct. Spaces} {\bf 2018}, Art. ID 9610257, 7 pp.



\bibitem {CuiW14} Cui, H. H.; Wang, F. H. Iterative methods for the split common fixed point problem in Hilbert spaces. {\it Fixed Point Theory Appl.} {\bf 2014}, 2014:78, 8 pp.

\bibitem {CuiW23} Cui, H. H.; Wang, F. H. The split common fixed point problem with multiple output sets for demicontractive mappings.
{\it Optimization} Published Online: 21 Feb 2023.

\bibitem {Dang} Dang, Y. Z.; Meng, F. W.; Sun, J. An iterative algorithm for split common fixed-point problem for demicontractive mappings. in {\it Optimization methods, theory and applications}, 85--94, Springer, Heidelberg, 2015.



\bibitem{Eslam16} Eslamian, M. General algorithms for split common fixed point problem of demicontractive mappings. {\it Optimization} {\bf 65} (2016), no. 2, 443--465.

\bibitem{Eslam17} Eslamian, M. Split common fixed point and common null point problem. {\it Math. Methods Appl. Sci.} {\bf 40} (2017), no. 18, 7410--7424.

\bibitem{Eslam17a} Eslamian, M.; Eskandani, G. Zamani; Raeisi, M. Split common null point and common fixed point problems between Banach spaces and Hilbert spaces. Mediterr. {\it J. Math.} {\bf 14} (2017), no. 3, Paper No. 119, 15 pp.

\bibitem{Fan21} Fan, Q.; Peng, J.; He, H. Weak and strong convergence theorems for the split common fixed point problem with demicontractive operators. {\it Optimization} {\bf 70} (2021), no. 5-6, 1409--1423.

\bibitem{Fan23} Fan, H. L.; Wang, C. Stability and convergence rate of Jungck-type iterations for a pair of strongly demicontractive mappings in Hilbert spaces. {\it Comput. Appl. Math.} {\bf 42} (2023), no. 1, Paper No. 33, 17 pp.


\bibitem{Gubin} Gubin, L. G.; Polyak, B. T.; Rajk, E. V. The method of projections for finding the common point of convex sets. {\it U.S.S.R. Comput.
Math. Math. Phys.} {\bf 7} (1967), No. 6, 1--24 (1970); translation from {\it Zh. Vychisl. Mat. Mat. Fiz.} {\bf 7}, 1211--1228 (1967).

\bibitem{Gupta} Gupta N.; Postolache M.; Nandal A.; Chugh, R. A cyclic iterative algorithm for multiple-sets split common fixed point problem of
demicontractive mappings without prior knowledge of operator norm. {\it Math.} {\bf 9} (2021), no. 4, 1--9, Article number 372.

\bibitem {Han18} Hanjing, A.; Suantai, S. Solving split equality common fixed point problem for infinite families of demicontractive mappings. 
 {\it Carpathian J. Math.} {\bf 34} (2018), no. 3, 321--331.

\bibitem {Han18a} Hanjing, A.; Suantai, S.  The split common fixed point problem for infinite families of demicontractive mappings. {\it Fixed Point Theory Appl.} {\bf 2018}, Paper No. 14, 21 p.

\bibitem {Han20} Hanjing, A.; Suantai, S. The split fixed point problem for demicontractive mappings and applications. {\it Fixed Point Theory}  {\bf 21} (2020), no. 2, 507--524.

\bibitem {Han20a} Hanjing, A.; Suantai, S.  Hybrid inertial accelerated algorithms for split fixed point problems of demicontractive mappings and equilibrium problems. {\it Numer. Algorithms}  {\bf 85} (2020), no. 3, 1051--1073.

\bibitem{Hanjing23} Hanjing, A.; Suantai, S.; Cho, Y. J. Hybrid inertial accelerated extragradient algorithms for split pseudomonotone equilibrium problems and fixed point problems of demicontractive mappings. {\it Filomat} {\bf 37} (2023), no. 5, 1607--1623.

\bibitem{He21a} He, H. M.; Fan, Q. W.; Chen, R. D. A new iterative construction for approximating solutions of a split common fixed point problem. {\it J. Math.} {\bf 2021}, Art. ID 6659079, 12 pp.

\bibitem{He16} He, H. M.; Liu, S. Y.; Chen, R. D. Strong convergence theorems for an implicit iterative algorithm for the split common fixed point problem. {\it J. Funct. Spaces} {\bf 2016}, Art. ID 4093524, 7 pp.

\bibitem{He16a} He, H. M.; Liu, S. Y.; Chen, R. D.; Wang, X. Y. Strong convergence results for the split common fixed point problem. {\it J. Nonlinear Sci. Appl.} {\bf 9} (2016), no. 9, 5332--5343.

\bibitem{He21} He, H. M.; Peng, J.; Fan, Q. W. An iterative viscosity approximation method for the split common fixed-point problem. {\it Optimization} {\bf 70} (2021), no. 5-6, 1261--1274.

\bibitem{HeDu14} He, Z. H.; Du, W.-S. On split common solution problems: new nonlinear feasible algorithms, strong convergence results and their applications. {\it Fixed Point Theory Appl.} {\bf 2014}, 2014:219, 16 pp.

\bibitem {Hicks} Hicks, T. L.; Kubicek, J. D. On the Mann iteration process in a Hilbert space. {\em J. Math. Anal. Appl.} {\bf 59} (1977), no. 3, 498--504.

\bibitem {Isio} Isiogugu, F.  O.; Pillay, P.; Baboolal, D. Approximation of a common element of the set of fixed points of multi-valued type-one demicontractive-type mappings and the set of solutions of an equilibrium problem in Hilbert spaces. {\em J. Nonlinear Convex Anal.} {\bf 17} (2016), no. 6, 1181--1197.

\bibitem {Jail} Jailoka, P.; Berinde, V.; Suantai, S. Strong convergence of Picard and Mann iterations for strongly demicontractive multi-valued mappings. {\it Carpathian J. Math.} {\bf 36} (2020), no. 2, 269--276.

\bibitem{Jail18} Jailoka, P.; Suantai, S. Split null point problems and fixed point problems for demicontractive multivalued mappings. {\it Mediterr. J. Math.} {\bf 15} (2018), no. 5, Paper No. 204, 19 pp.

\bibitem{Jail19} Jailoka, P.; Suantai, S. Split common fixed point and null point problems for demicontractive operators in Hilbert spaces. {\it Optim. Methods Softw.} {\bf 34} (2019), no. 2, 248--263.

\bibitem{Jail20} Jailoka, P.; Suantai, S. Viscosity approximation methods for split common fixed point problems without prior knowledge of the operator norm. {\it Filomat} {\bf 34} (2020), no. 3, 761--777.

\bibitem{Jail19a} Jailoka, P.; Suantai, S. The split common fixed point problem for multivalued demicontractive mappings and its applications. {\it Rev. R. Acad. Cienc. Exactas F\' is. Nat. Ser. A Mat. RACSAM} {\bf 113} (2019), no. 2, 689--706.

\bibitem{Jail21} Jailoka, P.; Suantai, S.  On split fixed point problems for multi-valued mappings and designing a self-adaptive method. {\it Results Math.} {\bf 76} (2021), no. 3, Paper No. 133, 21 pp.

\bibitem{Jail21} Jailoka, P.; Suantai, S.; Hanjing, A. A fast viscosity forward-backward algorithm for convex minimization problems with an application in image recovery. {\it Carpathian J. Math.} {\bf 37} (2021), no. 3, 449--461.

\bibitem{Jira} Jirakitpuwapat W.; Kumam P.; Cho Y. J.; Sitthithakerngkiet K. A general algorithm for the split common fixed point problem with its applications to signal processing. {\it Math.} {\bf 7} (2019), no. 31 Article number 226.


\bibitem {Kimura} Kimura, Y. Resolvents of equilibrium problems on a complete geodesic space with curvature bounded above. {\it Carpathian J. Math.} {\bf 37} (2021), no. 3, 463--476.

\bibitem {Kitkuan} Kitkuan, D.; Kumam, P.; Berinde, V.; Padcharoen, A. Adaptive algorithm for solving the SCFPP of demicontractive operators without a priori knowledge of operator norms. {\it An. \c Stiin\c t. Univ. "Ovidius" Constan\c ta Ser. Mat.} {\bf 27} (2019), no. 3, 153--175.

\bibitem {Lin} Lin, L.-J. Bilevel problems over split equality fixed point for finite families of countable nonlinear mappings. {\it J. Nonlinear Convex Anal.} {\bf 21} (2020), no. 1, 221--241.

\bibitem {Linh} Linh, H. M.; Reich, S.; Thong, D. V.; Dung, V. T.; Lan, N. P. Analysis of two variants of an inertial projection algorithm for finding the minimum-norm solutions of variational inequality and fixed point problems. {\it Numer. Algorithms} {\bf 89} (2022), no. 4, 1695--1721.


\bibitem {Majee} Majee, P.; Bai, S. N.; Padhye, S. Inertial Mann type algorithms for a finite collection of equilibrium problems and fixed point problem of demicontractive mappings. {\it J. Analysis.} Sep 2023 (Early Access)

\bibitem {Main08} Maing\' e, P.-E.  Convex minimization over the fixed point set of demicontractive mappings. {\em Positivity} {\bf 12} (2008), no. 2, 269--280.

\bibitem {Main08a} Maing\' e, P.-E. A hybrid extragradient-viscosity method for monotone operators and fixed point problems. {\em SIAM J. Control Optim.} {\bf 47} (2008), no. 3, 1499--1515.

\bibitem {Main08d} Maing\' e, P.-E.  Extension of the hybrid steepest descent method to a class of variational inequalities and fixed point problems with nonself-mappings. {\it Numer. Funct. Anal. Optim.} {\bf 29} (2008), no. 7-8, 820--834.

\bibitem {Main08b} Maing\' e, P.-E.  Regularized and inertial algorithms for common fixed points of nonlinear operators. {\it J. Math. Anal. Appl.} {\bf 344} (2008), no. 2, 876--887.

\bibitem {Main08c} Maing\' e, P.-E.  New approach to solving a system of variational inequalities and hierarchical problems. {\it J. Optim. Theory Appl.} {\bf 138} (2008), no. 3, 459--477.

\bibitem {Main11} Maing\' e, P.-E.; M\u aru\c ster, \c St.  Convergence in norm of modified Krasnoselski-Mann iterations for fixed points of demicontractive mappings. {\it Appl. Math. Comput.} {\bf 217} (2011), no. 24, 9864--9874.

\bibitem {Main08e} Maing\' e, P.-E.;  Moudafi, A. Coupling viscosity methods with the extragradient algorithm for solving equilibrium problems. {\it J. Nonlinear Convex Anal.} {\bf 9} (2008), no. 2, 283--294.

\bibitem {Marino}  Marino, G.; Xu, H.-K. Weak and strong convergence theorems for strict pseudo-contractions in Hilbert spaces. {\it J. Math. Anal. Appl.} {\bf 329} (2007), no. 1, 336--346.

\bibitem {MarMar05} M\u aru\c ster, L.; M\u aru\c ster, \c St. On convex feasibility problems. {\it Carpathian J. Math.} {\bf 21} (2005), no. 1-2, 83--87.

\bibitem {Mar73} M\u aru\c ster, \c St.  Sur le calcul des z\' eros d'un op\' erateur discontinu par it\' eration. {\it Canad. Math. Bull.} {\bf 16} (1973), 541--544.

\bibitem {Mar77} M\u aru\c ster, \c St. The solution by iteration of nonlinear equations in Hilbert spaces. {\it Proc. Amer. Math. Soc.} {\bf 63} (1977), no. 1, 69--73.

\bibitem {Mar03} M\u aru\c ster, \c St. On the projection methods for convex feasibility problems. {\it An. Univ. Timi\c soara Ser. Mat.-Inform.} {\bf 41} (2003), Special issue, 177--182.

\bibitem {Mar05} M\u aru\c ster, \c St.  Quasi-nonexpansivity and the convex feasibility problem. {\it An. \c Stiin\c t. Univ. Al. I. Cuza Ia\c si Inform. (N.S.)} {\bf 15} (2005), 47--56 (2006).

\bibitem {MarP08} M\u aru\c ster, \c St.; Popirlan, C. On the Mann-type iteration and the convex feasibility problem. {\it J. Comput. Appl. Math.} {\bf 212} (2008), no. 2, 390--396.

\bibitem {MarP09} M\u aru\c ster, \c St.; Popirlan, C. On the regularity condition in a convex feasibility problem. Nonlinear Anal. 70 (2009), no. 5, 1923--1928.

\bibitem {MarR15} M\u aru\c ster, \c St. Rus, I. A. Kannan contractions and strongly demicontractive mappings. {\it Creat. Math. Inform.} {\bf 24} (2015), no. 2, 171--180.


\bibitem{Meddahi} Meddahi, M.; Nachi, K.; Benahmed, B. A hybrid conjugate method for variational inequalities over fixed point sets of demicontractive multimaps. {\it Nonlinear Stud.} {\bf 27} (2020), no. 4, 975--989.

\bibitem {Mewomo} Mewomo, O. T.; Ogbuisi, F. U.; Okeke, C. C. On split equality minimization and fixed point problems. {\it Novi Sad J. Math.} {\bf 48}
(2018), no. 2, 21--39.

\bibitem{minjibir2022converg} Minjibir, M. S.; Salisu, S. Strong and $\Delta$-convergence theorems for a countable family of multivalued demicontractive maps in Hadamard spaces. {\it Nonlinear Funct. Anal. Appl.}, {\bf 27}(1) (2022), 45--58.

\bibitem {Moore} Moore, C., {\it Iterative approximation fixed points of demicontractive maps}, The Abdus Salam Intern. Centre for Theoretical Physics,Trieste, Italy, Scientific Report, IC/98/214, November, 1998.

\bibitem{Moudafi} Moudafi, A. The split common fixed-point problem for demicontractive mappings. {\it Inverse Problems} {\bf 26} (2010), no. 5, 055007, 6 pp.

\bibitem{Mongkol} Mongkolkeha, C.; Cho, Y. J.; Kumam, P. Convergence theorems for $k$-dimeicontactive mappings in Hilbert spaces. [[Convergence theorems for $k$-demicontractive mappings in Hilbert spaces]] {\it Math. Inequal. Appl.} {\bf 16} (2013), no. 4, 1065--1082.

\bibitem{Mouk} Mouktonglang, T.; Suparatulatorn, R. Inertial hybrid projection methods with selection techniques for split common fixed point problems in Hilbert spaces. {\it Politehn. Univ. Bucharest Sci. Bull. Ser. A Appl. Math. Phys.} {\bf 84} (2022), no. 2, 47--54.

\bibitem {Mouk23} Mouktonglang, Thanasak; Poochinapan, Kanyuta; Varnakovida, Pariwate; Suparatulatorn, Raweerote; Moonchai, Sompop. Convergence analysis of two parallel methods for common variational inclusion problems involving demicontractive mappings. {\it J. Math.} {\bf 2023}, Art. ID 1910411, 19 pp.


\bibitem{Ogbuisi20} Ogbuisi, F. U.; Isiogugu, F. O. A new iterative algorithm for pseudomonotone equilibrium problem and a finite family of demicontractive mappings. {\it Abstr. Appl. Anal.} {\bf 2020}, Art. ID 3183529, 11 pp.

\bibitem{Ogbuisi19} Ogbuisi F.O.; Mewomo O.T. Strong convergence result for solving split hierarchical variational inequality problem for demicon-
tractive mappings. {\it Adv. Nonlinear Var. Inequal.} {\bf 22} (2019), no. 1, 24--39.

\bibitem{Okeke20} Okeke, C. C.; Izuchukwu, C.; Mewomo, O. T. Strong convergence results for convex minimization and monotone variational
inclusion problems in Hilbert space. {\it Rend. Circ. Mat. Palermo (2)} {\bf 69} (2020), no. 2, 675--693.

\bibitem{Okeke22} Okeke, C. C.; Ugwunnadi, G. C.; Jolaoso, L. O. An extragradient inertial algorithm for solving split fixed-point problems of
demicontractive mappings, with equilibrium and variational inequality problems. {\it Demonstr. Math.} {\bf 55} (2022), no. 1, 506--527.

\bibitem {Osilike93} Osilike, M. O. Iterative method for nonlinear monotone-type operators in uniformly smooth Banach spaces. {\it J. Nigerian Math. Soc.} {\bf 12} (1993), 73--79.

\bibitem {Osilike00} Osilike, M. O. Strong and weak convergence of the Ishikawa iteration method for a class of nonlinear equations. {\it Bull. Korean Math. Soc.} {\bf 37} (2000), no. 1, 153--169.

\bibitem{Pad21} Padcharoen, A.; Kumam, P.; Cho, Y. J. Split common fixed point problems for demicontractive operators. {\it Numer. Algorithms} {\bf 82} (2019), no. 1, 297--320.

\bibitem{Pany23} Panyanak, B.; Khunpanuk, C.; Pholasa, N.; Pakkaranang, N. Dynamical inertial extragradient techniques for solving equilibrium and fixed-point problems in real Hilbert spaces. {\it J. Inequal. Appl.} {\bf 2023}, Paper No. 7, 36 pp.



\bibitem{Pet} Petryshyn, W. V. Iterative construction of fixed points of contractive type mappings in Banach spaces. 1968 {\it Numerical Analysis of Partial Differential Equations} (C.I.M.E. 2 Ciclo, Ispra, 1967) pp. 307--339 Edizioni Cremonese, Rome. (re-edited by Spinger in 2011).

\bibitem {Qin} Qin, L.-J.; Wang, G. Multiple-set split feasibility problems for a finite family of demicontractive mappings in Hilbert spaces. {\it Math. Inequal. Appl.} {\bf 16} (2013), no. 4, 1151--1157.

\bibitem {Rehman} Rehman, H. ur; Kumam, P.; Kumam, W.; Sombut, K. A new class of inertial algorithms with monotonic step sizes for solving fixed point and variational inequalities. {\it Math. Methods Appl. Sci.} {\bf 45} (2022), no. 16, 9061--9088.

\bibitem {Rehman23} Rehman, H. ur; Kumam, P.; Berinde, V. A family of monotonic iterative methods for solving $\rho$-demicontractive fixed point problems and variational inequalities involving pseudomonotone operators. {\it J. Nonlinear Convex Anal.} {\bf 24} (2023), no. 4, 905--924.

\bibitem {Saejung} Saejung, S.; Kraikaew, R. A unified algorithm for finding a fixed point of demicontractive mappings and its application to split common fixed point problem. {\it Rev. R. Acad. Cienc. Exactas F\' is. Nat. Ser. A Mat. RACSAM} {\bf 116} (2022), no. 4, Paper No. 190, 11 pp.

\bibitem {Salisu} Salisu, S.; Berinde, V.; Sriwongsa, S.; Kumam, P. Approximating fixed points of demicontractive mappings in metric spaces by geodesic averaged perturbation techniques. {\it AIMS Math.} {\bf 8} (2023), no. 12, 28582--28600.



\bibitem {Shehu15} Shehu, Y. New convergence theorems for split common fixed point problems in Hilbert spaces. {\it J. Nonlinear Convex Anal.} {\bf 16} (2015), no. 1, 167--181.

\bibitem {Shehu16} Shehu, Y.; Cholamjiak, P. Another look at the split common fixed point problem for demicontractive operators. {\it Rev. R. Acad. Cienc. Exactas F\' is. Nat. Ser. A Mat. RACSAM} {\bf 110} (2016), no. 1, 201--218.

\bibitem {Shehu16a} Shehu, Y.; Mewomo, O. T. Further investigation into split common fixed point problem for demicontractive operators. {\it Acta Math. Sin. (Engl. Ser.)} {\bf 32} (2016), no. 11, 1357--1376.

\bibitem{sow2021alg} Sow, T. M. M. Algorithm for computing a common solution of equilibrium and fixed point problems with set-valued demicontractive operators. {\it Khayyam J. Math.} {\bf 7} (2021), no. 1, 131--139.

\bibitem {Sow20b} Sow, T. M. M. General iterative algorithm for demicontractive-type mapping in real Hilbert spaces. {\it Creat. Math. Inform.} {\bf 29} (2020), no. 1, 91--99.

\bibitem {Sow20a} Sow, T. M. M. A new iterative algorithm for solving some nonlinear problems in Hilbert spaces. {\it J. Nonlinear Sci. Appl.} {\bf 13} (2020), no. 3, 119--132.

\bibitem {Sow21a} Sow, T. M. M. Nonlinear iterative algorithms for solving variational inequality problems over the set of common fixed point of one-parameter nonexpansive semigroup and demicontractive mappings. {\it Asian-Eur. J. Math.}  {\bf 14} (2021), no. 10, Paper No. 2150170, 18 pp.

\bibitem {Sow22} Sow, T. M. M. A modified forward-backward splitting method for sum of monotone operators and demicontractive mappings. {\it Casp. J. Math. Sci.} {\bf 11} (2022), no. 1, 229--241.

\bibitem {Su} Su, H. Y.; Zhao, J. Self-adaptive iterative algorithms for solving multiple-set split equality common fixed-point problems of demicontractive operators. {\it J. Nonlinear Funct. Anal.} {\bf 2018} (2018), Article ID 47, 1--17.

\bibitem {Sua20} Suantai, S.; Jailoka, P. A self-adaptive algorithm for split null point problems and fixed point problems for demicontractive multivalued mappings. {\it Acta Appl. Math.} {\bf 170} (2020), 883--901. 


\bibitem {Sua17} Suantai, S.; Phuengrattana, W. A hybrid shrinking projection method for common fixed points of a finite family of demicontractive mappings with variational inequality problems. {\it Banach J. Math. Anal.} {\bf 11} (2017), no. 3, 661--675.

\bibitem {Sua21} Suantai, S.;  Sarnmeta, P.; Chumpungam, D.; Inthakon, W. Split common fixed point problems for multi-valued demicontractive mappings in Hilbert spaces. {\it J. Nonlinear Convex Anal.} {\bf 22} (2021), no. 12, 2623--2637.

\bibitem{Supa20a} Suparatulatorn, R. Weak convergence theorem of generalized self-adaptive algorithms for solving split common fixed point problems. {\it Politehn. Univ. Bucharest Sci. Bull. Ser. A Appl. Math. Phys.} {\bf 82} (2020), no. 2, 67--74.

\bibitem{Supa20} Suparatulatorn, R.; Cholamjiak, P.; Suantai, S. Self-adaptive algorithms with inertial effects for solving the split problem of the demicontractive operators. {\it Rev. R. Acad. Cienc. Exactas F\' is. Nat. Ser. A Mat. RACSAM} {\it 114} (2020), no. 1, Paper No. 40, 16 pp.

\bibitem{Supa19} Suparatulatorn, R.; Charoensawan, P.; Poochinapan, K. Inertial self-adaptive algorithm for solving split feasible problems with applications to image restoration. {\it Math. Methods Appl. Sci.} {\bf 42} (2019), no. 18, 7268--7284.


\bibitem{Supa21} Suparatulatorn, R.; Charoensawan, P.; Poochinapan, K.; Dangskul, S. An algorithm for the split feasible problem and image restoration. {\it Rev. R. Acad. Cienc. Exactas F\' is. Nat. Ser. A Mat. RACSAM} {\bf 115} (2021), no. 1, Paper No. 12, 18 pp.

\bibitem{Supa20b} Suparatulatorn, R.; Khemphet, A.; Charoensawan, P.; Suantai, S.; Phudolsitthiphat, N. Generalized self-adaptive algorithm for solving split common fixed point problem and its application to image restoration problem. {\it Int. J. Comput. Math.} {\bf 97} (2020), no. 7, 1431--1443.

\bibitem{Supa19a} Suparatulatorn, R.; Suantai, S.; Phudolsitthiphat, N. Reckoning solution of split common fixed point problems by using inertial self-adaptive algorithms. {\it Rev. R. Acad. Cienc. Exactas F\' is. Nat. Ser. A Mat. RACSAM} {\bf 113} (2019), no. 4, 3101--3114.

\bibitem{Tan22} Tan, B.; Liu, L.; Qin, X. Strong convergence of inertial extragradient algorithms for solving variational inequalities and fixed point problems. {\it Fixed Point Theory} {\bf 23} (2022), no. 2, 707--727.

\bibitem{Tan22a} Tan, B.; Zhou, Z.; Li, S. X. Viscosity-type inertial extragradient algorithms for solving variational inequality problems and fixed point problems. {\it J. Appl. Math. Comput.} {\bf 68} (2022), no. 2, 1387--1411.

\bibitem{Tang12} Tang, Y.-C.; Peng, J.-G.; Liu, L.-W.  A cyclic algorithm for the split common fixed point problem of demicontractive mappings in Hilbert spaces. {\it Math. Model. Anal.} {\bf 17} (2012), no. 4, 457--466.

\bibitem{Tang14} Tang, Y.-C.; Peng, J.-G.; Liu, L.-W. A cyclic and simultaneous iterative algorithm for the multiple split common fixed point problem of demicontractive mappings. {\it Bull. Korean Math. Soc.} {\bf 51} (2014), no. 5, 1527--1538.

\bibitem {Thong22} Thong, D. V.; Dung, V. T.; Long, L. V. Inertial projection methods for finding a minimum-norm solution of pseudomonotone variational inequality and fixed-point problems. {\bf Comput. Appl. Math.} {\bf 41} (2022), no. 6, Paper No. 254, 25 pp

\bibitem {Thong18} Thong, D. V.; Hieu, D. V. Modified subgradient extragradient algorithms for variational inequality problems and fixed point problems.  {\it Optimization} {\bf 67} (2018), no. 1, 83--102.

\bibitem {Thong18a} Thong, D. V.; Hieu, D. V.  A new approximation method for finding common fixed points of families of demicontractive operators and applications. {\it J. Fixed Point Theory Appl.} {\bf 20} (2018), no. 2, Paper No. 73, 27 pp.

\bibitem {Thong23} Thong, D. V.; Liu, L.-L.; Dong, Q.-L.; Van Long, L.; Tuan, P. A. Fast relaxed inertial Tseng's method-based algorithm for solving variational inequality and fixed point problems in Hilbert spaces. {\it J. Comput. Appl. Math.} {\bf 418} (2023), Paper No. 114739, 22 pp.

\bibitem {Tic} \c Tical\u a, C. Approximating fixed points of demicontractive mappings by iterative methods defined as admissible perturbations. {\em Creat. Math. Inform.} {\bf 25} (2016), no. 1, 121--126.

\bibitem{Vuong} Vuong, P. T.; Strodiot, J. J.; Nguyen, V. H. On extragradient-viscosity methods for solving equilibrium and fixed point problems in a Hilbert space.  {\it Optimization} {\bf 64} (2015), no. 2, 429--451.

\bibitem{Wang21} Wang, A.; Zhao, J. Self-adaptive iterative algorithms for the split common fixed point problem with demicontractive operators. {\it J.
Nonlinear Var. Anal.} {\bf 5} (2021), no. 4, 573--587.

\bibitem{Wang-Fan} Wang, C.; Fan, H. L. Weak stability of Ishikawa iterations for strongly demicontractive mappings in Hilbert spaces. {\it Filomat} {\bf 36} (2022), no. 14, 4869--4876.

\bibitem{Wang-multi} Wang, C. S.; Ceng, L. C.; Li, B.; Cao, S.-L.; HU, H. -Y.; Liang, Y. S. Modified Inertial-Type Subgradient Extragradient Methods for Variational Inequalities and Fixed Points of Finite Bregman Relatively Nonexpansive and Demicontractive Mappings. {\it Axioms} {\bf 2023}, 12(9), 832; https://doi.org/10.3390/axioms12090832

\bibitem{Wang22} Wang, F. H. The split feasibility problem with multiple output sets for demicontractive mappings. {\it J. Optim. Theory Appl.} {\bf 195} (2022), no. 3, 837--853.

\bibitem{Wang13} Wang, F. H.; Cui, H. H. Convergence of a cyclic algorithm for the split common fixed point problem without continuity assumption. {\it Math. Model. Anal.} {\bf 18} (2013), no. 4, 537--542.

\bibitem{Wang19} Wang, J. Y.; Fang, X. L. A strong convergence theorem for the split common fixed-point problem of demicontractive mappings.
{\it Appl. Set-Valued Anal. Optim.} {\bf 1} (2019), No. 2, 105--112.

\bibitem{Wang20a} Wang, Y. Q.; Chen, J.; Pitea, A. The split equality fixed point problem of demicontractive operators with numerical example and
application. {\it Symmetry-Basel} {\bf 12} (2020), no. 6, Article number 902.

\bibitem{Wang20} Wang, Y. Q.; Fang, X. L.; Kim, T.-H. A new algorithm for common fixed-point problems of a finite family of asymptotically demicontractive operators and its applications. {\it J. Nonlinear Convex Anal.} {\bf 21} (2020), no. 9, 1875--1887.

\bibitem{Wang17b} Wang, Y. Q.; Fang, X. L. Viscosity approximation methods for the multiple-set split equality common fixed-point problems of demicontractive mappings. {\it J. Nonlinear Sci. Appl.} {\bf 10} (2017), no. 8, 4254--4268.

\bibitem{Wang21} Wang, Y. Q.; Fang, X. L.; Guan, J.-L.; Kim, T.-H. On split null point and common fixed point problems for multivalued demicontractive mappings. {\em Optimization} {\bf 70} (2021), no. 5-6, 1121--1140.

\bibitem{Wang17c} Wang, Y. Q.; Kim, T.-H.;  Fang, X. L.; He, H. M. The split common fixed-point problem for demicontractive mappings and quasi-nonexpansive mappings. {\it J. Nonlinear Sci. Appl.} {\bf 10} (2017), no. 6, 2976--2985.

\bibitem{Wang17} Wang, Y. Q.; Kim, T.-H.; Chen, R. D.; Fang, X. L. The multiple-set split equality common fixed point problems for demicontractive mappings without prior knowledge of operator norms. {\it J. Nonlinear Convex Anal.} {\bf 18} (2017), no. 10, 1849--1865.

\bibitem{Wang17a} Wang, Y. Q.; Kim, T.-H.;  Fang, X. L.  Weak and strong convergence theorems for the multiple-set split equality common fixed-point problems of demicontractive mappings. {\it J. Funct. Spaces} {\bf 2017}, Art. ID 5306802, 11 pp.

\bibitem{Wang18} Wang, Y. Q.; Liu, W.; Song, Y. L.; Fang, X. L. Mixed iterative algorithms for the multiple-set split equality common fixed-point problem of demicontractive mappings. {\it J. Nonlinear Convex Anal.} {\bf 19} (2018), no. 11, 1921--1932.

\bibitem{Wang18a} Wang, J. Q.; Wang, Y. Strong convergence of a cyclic iterative algorithm for split common fixed-point problems of demicontractive
mappings. {\it J. Nonlinear Var. Anal.} {\bf 2} (2018), no. 3, 295--303.

\bibitem{Yang} Yang, L.; Zhao, F.; Kim, J. K. The split common fixed point problem for demicontractive mappings in Banach spaces. {\it J. Comput. Anal. Appl.} {\bf 22} (2017), no. 5, 858--863.

\bibitem{Xiao20} Xiao, J. F.; Huang, L.; Wang, Y. Q. Strong convergence of modified inertial Halpern simultaneous algorithms for a finite family of
demicontractive mappings. {\it Appl. Set-Valued Anal. Optim.} {\bf 2} (2020), No. 3, pp. 317--327.

\bibitem{Xiao22} Xiao, J. F.; Wang, Y. Q. A viscosity method with inertial effects for split common fixed point problems of demicontractive mappings.
{\it J. Nonlinear Funct. Anal.} {\bf 2022} (2022), Article number 17.

\bibitem{Xu} Xu, H. Y.; Lan, H. Y. Novel extended Halpern-type convergence algorithms for the split common fixed point problem involving a-demicontractive operators. {\it Mathematics} {\bf 2023}, 11(18), 3871; https://doi.org/10.3390/math11183871.

\bibitem{Yao18a} Yao, Y. H.; Leng, L. M.; Liou, Y.-C. Strong convergence of an iteration for the split common fixed points of demicontractive operators. {\it J. Nonlinear Convex Anal.} {\bf 19} (2018), no. 2, 197--205.


\bibitem{Yao18} Yao, Y. H.; Liou, Y.-C.;  Wu, Y.-J. An extragradient method for mixed equilibrium problems and fixed point problems. {\it Fixed Point Theory Appl.} {\bf 2009}, Art. ID 632819, 15 pp.

\bibitem{Yao18b} Yao, Y. H.; Liou, Y.-C.; Postolache, M. Self-adaptive algorithms for the split problem of the demicontractive operators. {\it Optimization} {\bf 67} (2018), no. 9, 1309--1319.

\bibitem{Yao18c} Yao, Y. H.; Qin, X. L.; Yao, J.-C. Self-adaptive step-sizes choice for split common fixed point problems. {\it J. Nonlinear Convex Anal.} {\bf 19} (2018), no. 11, 1959--1969.

\bibitem{Yao20} Yao, Y. H.; Yao, J.-C.; Liou, Y.-C.; Postolache, M. Iterative algorithms for split common fixed points of demicontractive operators without priori knowledge of operator norms. {\it Carpathian J. Math.} {\bf 34} (2018), no. 3, 459--466.

\bibitem{Ying} Ying, Y.; Huang, L.; Zhang, Y. Q. Viscosity approximation of a modified inertial simultaneous algorithm for a finite family of
demicontractive mappings. {\it J. Nonlinear Funct. Anal.} {\bf 2023} (2023) Article number 7.

\bibitem{Yu21} Yu, Y. L. Analysis of algorithms for solving variational inclusions and split fixed point problems. {\it J. Nonlinear Convex Anal.} {\bf 22} (2021), no. 1, 87--96.

\bibitem{Yu12} Yu, Y. R.; Sheng, D. L. On the strong convergence of an algorithm about firmly pseudo-demicontractive mappings for the split common fixed-point problem. {\it J. Appl. Math.} {\bf 2012}, Art. ID 256930, 9 pp.

\bibitem{Zhao} Zhao, Y. F.; Yao,Y. H. Weak convergence of a new iterate for solving split fixed point problems. {\it Politehn.Univ. Bucharest Sci.
Bull. Ser. A Appl. Math. Phys.} {\bf 85} (2023), no. 2, 43--50.

\bibitem{Zhang} Zhang, C. J.; Li, Y.; Wang, Y. H. On solving split generalized equilibrium problems with trifunctions and fixed point problems of demicontractive multi-valued mappings. {\em J. Nonlinear Convex Anal.} {\bf 21} (2020), no. 9, 2027--2042.

\bibitem{Zheng} Zheng, X. X.; Yao, Y.H.; Liou, Y.-C.; Leng, L. M. Fixed point algorithms for the split problem of demicontractive operators. {\it J. Nonlinear Sci. Appl.} {\bf 10} (2017), no. 3, 1263--1269.

\bibitem {Zhou08} Zhou, H. Convergence theorems of fixed points for $\kappa$-strict pseudo-contractions in Hilbert spaces. {\it Nonlinear Anal.} {\bf 69} (2008), no. 2, 456--462.


\bibitem {Zhou21} Zhou, Z.; Tan, B.; Li, S. X. An accelerated hybrid projection method with a self-adaptive step-size sequence for solving split common fixed point problems. {\it Math. Methods Appl. Sci.} {\bf 44} (2021), no. 8, 7294--7303.

\bibitem {Zhu} Zhu, L. -J.; Yao, Y. H. Algorithms for approximating solutions of split variational inclusion and fixed-point problems. {\it Math.} {\bf 11}
(2023), no. 3, Article number 641.

\bibitem {Zhu24} Zhu, L. -J.; Yao, J. C.; Yao, Y. H. Approximating solutions of a split fixed point problem of demicontractive operators. {\it Carpathian J. Math.} {\bf 40} (2024), No. 1; 195--206 DOI10.37193/CJM.2024.01.14

\bibitem {Zong} Zong, C. X.; Tang, Y. C. Iterative methods for solving the split common fixed point problem of demicontractive mappings in Hilbert spaces. {\it J. Nonlinear Sci. Appl.} {\bf 11} (2018), no. 8, 960--970.

\end{thebibliography}
\end{document}